\documentclass{ifacconf}

\usepackage{natbib} 
\usepackage[centertags]{amsmath}
\usepackage{graphics}
\usepackage{amscd}
\usepackage{amsfonts}
\usepackage{amssymb}
\usepackage{helvet}
\usepackage{rotating}
\usepackage{epsfig}
\usepackage[latin1]{inputenc}
\usepackage{natbib}
\usepackage[english,german,american]{babel}
\usepackage{algorithm}

\usepackage{cancel}

\newtheorem{theorem}{Theorem}[section]

\newtheorem{proposition}[theorem]{Proposition}

\newtheorem{remark}[theorem]{Remark}
\newtheorem{example}[theorem]{Example}
\newtheorem{corollary}[theorem]{Corollary}

\def\K{{\mathcal{K}}}

\def\R{\mathbb{R}}
\def\N{\mathbb{N}}

\def\cK{\mathcal{K}}
\def\cL{\mathcal{L}}
\def\cKL{\mathcal{KL}}

\def\x{x}
\def\X{\mathbb{R}^n}
\def\bX{\mathbb{X}}
\def\u{u}
\def\U{\mathbb{R}^m}
\def\bU{\mathbb{U}}
\def\cU{\mathcal{U}}
\def\xu{\x_\u}
\def\xmu{\x_\mu}
\def\xs{\x^\star}
\def\xh{\hat{\x}}
\def\us{\u^\star}
\def\l{\ell}

\usepackage{color}
\definecolor{karl}{rgb}{.9,.1,.4} 
\definecolor{juergen}{rgb}{.1,.4,.9} 
\definecolor{philipp}{rgb}{.4,.9,.1} 
\definecolor{karlfrage}{rgb}{.2,.9,.9}

\begin{document}

% ------------------------------------------------------------------------
\begin{frontmatter}
%\title{Stability of Predictive Feedbacks with Short Optimization Horizon}
%\title{Stability of Prediction Based Algorithms despite Shortened Optimization Horizons}
\title{Predictive Control Algorithms: Stability despite Shortened Optimization Horizons}
%\title{Stability of Optimization based Algorithms with Time Varying Control Horizons}
\author[Second]{Philipp Braun} 
\author[First]{J\"{u}rgen Pannek} 
\author[Second]{Karl Worthmann} 

\address[Second]{University of Bayreuth, 95440 Bayreuth, Germany}%\newline (e-mail: \{karl.worthmann, philipp.braun\}@uni-bayreuth.de)}
\address[First]{University of the Federal Armed Forces, 85577 Munich, Germany}%\newline (e-mail: juergen.pannek@unibw.de)}


% ------------------------------------------------------------------------
\begin{abstract}
	: The stability analysis of model predictive control schemes without terminal constraints and/or costs has attracted considerable attention during the last years. We pursue a recently proposed approach which can be used to determine a suitable optimization horizon length for nonlinear control systems governed by ordinary differential equations. In this context, we firstly show how the essential growth assumption involved in this methodology can be derived and demonstrate this technique by means of a numerical example. Secondly, inspired by corresponding results, we develop an algorithm which allows to reduce the required optimization horizon length while maintaining asymptotic stability or a desired performance bound. Last, this basic algorithm is further elaborated in order to enhance its robustness properties.
\end{abstract}
% ------------------------------------------------------------------------

\end{frontmatter}
% ------------------------------------------------------------------------

% ------------------------------------------------------------------------
\section{Introduction}
% ------------------------------------------------------------------------

Within the last decades, model predictive control (MPC) has grown mature for both linear and nonlinear systems, see, e.g., \cite{CamachoBordons2004} or \cite{RawlingsMayne2009}. % or \cite{GruenePannek2011}. 
Although analytically and numerically challenging, the method itself is attractive due to its simplicity: In a first step, a new measurement of the current system state is obtained which is thereafter used to compute an optimal control over a finite optimization horizon. In the third and last step, a portion of this control is applied to the process and the entire problem is shifted forward in time rendering the scheme to be iteratively applicable.

Stability of the MPC closed loop can be shown by imposing endpoint constraints, Lyapunov type terminal costs or terminal regions, cf. \cite{KeerthiGilbert1988} and \cite{ChenAllgoewer1998}. Here, we study MPC schemes without these ingredients for which stability and, in addition, bounds on the required horizon length can be deduced, both for linear and nonlinear systems, cf.~\cite{PrimbsNevistic2000} and \cite{TunaMessinaTeel2006}. We follow the recent approach from \cite{RebleAllgoewer2011} extending \cite{Gruene2009,GruenePannekSeehaferWorthmann2010} %\cite{GruenePannekSeehaferWorthmann2010} 
to continuous time systems which not only guarantees stability but also reveals an estimate on the degree of suboptimality with respect to the optimal controller on an infinite horizon.

In this work, we show how the essential assumption needed to apply the methodology proposed in \cite{RebleAllgoewer2011} can be practically verified. Then, based on observations drawn from numerical computations, implementable MPC algorithms with variable control horizons are developed which allow for smaller optimization horizons while maintaining stability or a desired performance bound. To overcome the lack of robustness implied by prolonging the control horizon and, thus, staying in open loop for longer time intervals, conditions are presented which ensure that the control loop can be closed more often. Similar ideas were introduced in \cite{PannekWorthmann2011a} for a discrete time setting. Last, the computational effort is further reduced by introducing slack which allows to violate our main stability condition --- a relaxed Lyapunov inequality --- temporarily, cf.~\cite{Giselsson2010}.

The paper is organized as follows: In Section \ref{Section:Preliminaries} the problem formulation is given. In the ensuing Section \ref{Section:Stability and Performance Bounds}, we summarize stability results from \cite{RebleAllgoewer2011} and propose a technique to verify the key assumption which is illustrated by an example. Thereafter, we present algorithms which allow for shortening the optimization horizon by using time varying control horizons. Before drawing conclusions, Section \ref{Section:Slack Stability} contains results on how stability may be guaranteed by using weaker stability conditions.

% ------------------------------------------------------------------------
\section{Setup and Preliminaries}
\label{Section:Preliminaries}
% ------------------------------------------------------------------------

Let $\N$ and $\R$ denote the set of natural and real numbers respectively and $\| \cdot \|$ the Euclidean norm on $\R^d$, $d \in \N$. A continuous function $\eta: \R_{\geq 0} \rightarrow \R_{\geq 0}$ is called class $\cK_\infty$-function if it satisfies $\eta(0) = 0$, is strictly increasing and unbounded. A continuous function $\beta: \R_{\geq 0} \times \R_{\geq 0} \rightarrow \R_{\geq 0}$ is said to be of class $\cKL$ if for each $r > 0$ we have that $\lim_{t \rightarrow \infty} \beta(r, t) = 0$ holds, and for each $t \geq 0$ the condition $\beta(\cdot, t) \in \cK_\infty$ is satisfied.

Within this work we consider nonlinear time invariant control systems
\begin{align}
	\label{Preliminaries:eq:control system}
	\dot{\x}(t) = f(\x(t),\u(t))
\end{align}
where $\x(t) \in \R^n$ and $\u(t) \in \R^m$ denote the state and control at time $t \geq 0$. Constraints can be included via suitable subsets $\bX \subset \R^n$ and $\bU \subset \U$ of the state and control space, respectively. We denote a state trajectory which emanates from the initial state $\x_{0}$ and is subject to the control function $\u:\R_{\geq 0} \rightarrow \R^m$ by $\xu(t) = \xu(t;\x_{0})$. In the presence of constraints, a control function $\u$ is called admissible for $\x$ on the interval $[0,T)$ if the corresponding solution $\xu(\cdot; \x)$ exists and satisfies
\begin{align}
	\xu(t;\x) \in \bX, t \in [0,T], \quad \; \text{and} \; \quad \u(t) \in \bU, t \in [0,T).
\end{align}

The set of these admissible control functions is denoted by $\cU_{\x}([0,T))$. For an infinite time interval, $\u:\R_{\geq 0}\rightarrow \U$ is called admissible for $\x$ if $\u|_{[0,T)} \in \cU_{\x}([0,T))$ holds for each $T>0$ and the respective set is denoted by $\cU_{\x}([0,\infty))$. 

%For system \eqref{Preliminaries:eq:control system} we assume an equilibrium $\xs \in \bX$ to exist, i.e. there exists a control $\us \in \bU$ such that $f(\xs, \us) = 0$ holds. Our goal is to design a feedback control law $\mu: \R^n \rightarrow \R^m$ such that the resulting closed loop is asymptotically stable with respect to $\xs$, i.e. there exists $\beta \in \mathcal{KL}$ such that $\| \xmu(t;\x_0) - \xs \| \leq \beta(\| \x_0 - \xs\|,t)$, $t \geq 0$, holds for all $\x_0 \in \bX$ where $\xmu(\cdot;\x_0)$ denotes the closed loop trajectory induced by $\mu$. The stabilization task is to be accomplished in an optimal fashion which is measured by a cost functional. To this end, we introduce the continuous running cost $\l: \X \times \U \to \R_{\geq 0}$ satisfying
For system \eqref{Preliminaries:eq:control system} we assume an equilibrium $(\xs,\us) \in \bX \times \bU$ to exist, i.e. $f(\xs, \us) = 0$ holds. Our goal is to design a feedback control law $\mu: \R^n \rightarrow \R^m$ such that the resulting closed loop is asymptotically stable with respect to $\xs$, i.e. there exists $\beta \in \mathcal{KL}$ such that $\| \xmu(t;\x_0) - \xs \| \leq \beta(\| \x_0 - \xs\|,t)$, $t \geq 0$, holds for all $\x_0 \in \bX$ where $\xmu(\cdot;\x_0)$ denotes the closed loop trajectory induced by $\mu$. The stabilization task is to be accomplished in an optimal fashion which is measured by a cost functional. To this end, we introduce the continuous running cost $\l: \X \times \U \to \R_{\geq 0}$ satisfying
\begin{align*}
	\l(\xs, \us) = 0 \quad\text{and}\quad \inf_{u \in \bU} \l(\x, \u) > 0\ \forall\, \x \neq \xs.%\geq \underline{\eta}( \| \x - \xs \| )
\end{align*}
Then, for a given state $\x \in \bX$, the cost of an admissible control $\u \in \cU_\x([0,\infty))$ is
\begin{align*}
	J_\infty(\x, \u) := \int_0^\infty \l(\xu(t; \x), \u(t))\, dt.
\end{align*}
The computation of a corresponding minimizer is, in general, computationally hard due to the curse of dimensionality, cf. \cite{BardiCapuzzo1997}. Hence, we use model predictive control (MPC) to approximately solve this task. The central idea of MPC is to truncate the infinite horizon, i.e. to compute a minimizer $\us \in \cU_{\x}([0,T))$ of the cost functional
\begin{align}
	\label{Preliminaries:eq:finite horizon cost functional}
	J_T(\x, \u) := \int_0^T \l(\xu(t; \x), \u(t))\, dt,
\end{align}
which can be done efficiently using discretization methods and nonlinear optimization algorithms, see, e.g., \cite{Maciejowski2002} or \cite[Chapter 10]{GruenePannek2011}. % or \cite{NocedalWright2006}.
Furthermore, we define the corresponding optimal value function $V_T(\x) := \inf_{\u \in \cU_{\x}([0, T))} J_T(\x, \u)$, $T \in \R_{\geq 0} \cup \{ \infty \}$.

To obtain an infinite horizon control, only the first portion of the computed minimizer is applied, i.e. we define the feedback law via
\begin{align}
	\label{Preliminaries:eq:feedback}
	\mu_{T, \delta}(t; x) := \us(t), \quad t \in [0, \delta),
\end{align}
for the so called control horizon $\delta \in (0, T)$. Last, the optimal control problem is shifted forward in time which renders Algorithm \ref{Preliminaries:alg:mpc} to be iteratively applicable.

\begin{algorithm}
	\caption{MPC}\label{Preliminaries:alg:mpc}
	Given: $T > \delta > 0$
	\begin{enumerate}
		\item Measure the current state $\hat{\x}$
		\item Compute a minimizer $\us \in \cU_{\hat{\x}}([0, T))$ of \eqref{Preliminaries:eq:finite horizon cost functional} and define the MPC feedback law via \eqref{Preliminaries:eq:feedback}
		\item Implement $\mu_{T, \delta}(t; \hat{x})|_{t \in [0, \delta)}$, shift the horizon forward in time by $\delta$ and goto (1)
	\end{enumerate}
\end{algorithm}

The closed loop state trajectory emanating from the initial state $\x_0$ subject to the MPC feedback law $\mu_{T, \delta}$ from Algorithm \ref{Preliminaries:alg:mpc} is denoted by $\x_{\mu_{T, \delta}}(\cdot; \x_0)$. Furthermore, $u_{\text{MPC}}^{\mu_{T,\delta}}: \R_{\geq 0} \rightarrow \R^m$ denotes the control function obtained by concatenating the applied pieces of control functions, i.e.
\begin{equation*}
	\u_{\text{MPC}}^{\mu_{T,\delta}}(t;\x_0) = \mu_{T,\delta}(t-\lfloor t / \delta \rfloor \delta, \x_{\mu_{T,\delta}}(\lfloor t / \delta \rfloor \delta;\x_0)). 
	%\mu_{T,\delta}(t;\x_{\mu_{T,\delta}}(0; \x_0))|_{t \in [0,\delta)},\mu_{T,\delta}(t;\x_{\mu_{T,\delta}}(\delta; \x_0))|_{t \in [0,\delta)},\ldots.
\end{equation*}
The resulting MPC closed loop cost are given by
\begin{align*}
	V_\infty^{\mu_{T, \delta}}(\x) := \int_0^\infty \l(\x_{\mu_{T, \delta}}(t; \x), u_{\text{MPC}}^{\mu_{T, \delta}}(t))\, dt.
\end{align*}

Note that we tacitly assume that Problem \eqref{Preliminaries:eq:finite horizon cost functional} is solvable for all $x_0 \in \bX$ and the minimum is attained in each step of Algorithm \ref{Preliminaries:alg:mpc}. For a detailed discussion of feasibility issues we refer to \cite[Chapter 8]{GruenePannek2011}.

% ------------------------------------------------------------------------
\section{Stability and Performance Bounds}
\label{Section:Stability and Performance Bounds}
% ------------------------------------------------------------------------

Due to the truncation of the infinite horizon, stability and optimality properties of the optimal control may be lost. Yet, stability can be shown if the optimization horizon is sufficiently long, cf.~\cite{AlamirBornard1995, JadbabaieHauser2005}. Additionally, an optimization horizon length $T$ can be determined for which both asymptotic stability as well as a performance bound on the MPC closed loop in comparison to the infinite horizon control law hold.%, cf.~ \cite{RebleAllgoewer2011}. %To state this result, we introduce the optimal value function $V_T(\x) := \inf_{\u \in \cU_{\x}([0, T))} J_T(\x, \u)$ for $T \in \R_{\geq 0} \cup \{ \infty \}$ and the MPC closed loop cost 
% \begin{align*}
% 	V_\infty^{\mu_{T, \delta}}(\x) := \int_0^\infty \l(\x_{\mu_{T, \delta}}(t; \x), u_{\text{MPC}}^{\mu_{T, \delta}}(t))\, dt.
% \end{align*}
\begin{theorem}\label{Stability and Performance Bounds:thm:stability}
	Suppose a control horizon $\delta > 0$ and a monotone bounded function $B: \R_{\geq 0} \to \R_{\geq 0}$ satisfying
	\begin{align}
		\label{Stability and Performance Bounds:thm:stability:eq1}
		V_t(\x) \leq B(t) \inf_{u \in \bU} \l(\x, \u) =: B(t) \l^\star(\x) ,\quad t \geq 0,
	\end{align}
	for all $\x \in \bX$ to be given. If $T > \delta$ is chosen such that $\alpha_{T, \delta} > 0$ holds for
	\begin{align}
		\label{Stability and Performance Bounds:thm:stability:eq2}
		\hspace*{-0.1cm}\alpha_{T, \delta} := 1 - \frac{e^{-\int_\delta^T B(t)^{-1} dt} e^{-\int_{T-\delta}^T B(t)^{-1} dt}}{\left[ 1- e^{- \int_\delta^T B(t)^{-1} dt} \right] \left[ 1 - e^{-\int_{T-\delta}^T B(t)^{-1} dt} \right]},
		%\label{\hspace*{-0.1cm}\alpha_{T, \delta} := 1 - \frac{e^{-\int_\delta^T B(t)^{-1} dt}}{\left[ 1- e^{- \int_\delta^T B(t)^{-1} dt} \right] \left[ e^{\int_{T-\delta}^T B(t)^{-1} dt} -1 \right]},
	\end{align}
	then the relaxed Lyapunov inequality
	\begin{align}
		\label{Stability and Performance Bounds:thm:stability:eq3}
		V_T(x) \hspace*{-0.05cm} - \hspace*{-0.05cm} V_T(\x_{\mu_{T, \delta}}\hspace*{-0.05cm}(\delta; \x)) \geq \alpha_{T, \delta} \int_0^\delta \hspace*{-0.1cm} \l(\x_{\mu_{T, \delta}}\hspace*{-0.05cm}(t; \x), \mu_{T, \delta}(t; \x)) dt
	\end{align}
	as well as the performance estimate
	\begin{align}
		\label{Stability and Performance Bounds:thm:stability:eq4}
		V_\infty^{\mu_{T, \delta}}(\x) \leq \alpha_{T, \delta}^{-1} V_\infty(\x)
	\end{align}		
	are satisfied for all $\x \in \bX$. If, additionally, there exist $\cK_\infty$ functions $\underline{\eta}$, $\overline{\eta}$ such that $\underline{\eta}(x) \leq \l^\star(x) \leq \overline{\eta}(x)$ hold for all $\x \in \bX$, then the MPC closed loop is asymptotically stable for horizon length $T$.
\end{theorem}
A detailed proof of Theorem \ref{Stability and Performance Bounds:thm:stability} is given in \cite{RebleAllgoewer2011}. We point out that for a given control horizon $\delta > 0$ and a desired performance specification $\overline{\alpha} \in [0,1)$ on the MPC closed loop, there always exists an optimization horizon $T > \delta$ such that $\alpha_{T,\delta} > \overline{\alpha}$ is satisfied, cf.~\cite{Worthmann2012nmpc}. Interpreting \eqref{Stability and Performance Bounds:thm:stability:eq4} exemplarily, the choice $\overline{\alpha} = 0$ corresponds to asymptotic stability of the MPC closed loop whereas $\overline{\alpha} = 0.5$ limits the cost of the MPC control to double the cost of the infinite horizon control.

The crucial assumption which has to be verified in order to apply Theorem \ref{Stability and Performance Bounds:thm:stability} is the growth condition \eqref{Stability and Performance Bounds:thm:stability:eq1}. Here, we demonstrate this by means of the following example:
\begin{example}\label{Stability and Performance:ex:generator}
	Consider the system dynamics of a synchronous generator given by
	\begin{align*}
		\dot{x}_1(t) & = x_2(t) \\
		\dot{x}_2(t) & = - b_1 x_3(t) \sin(x_1(t)) - b_2 x_2(t) + P \\
		\dot{x}_3(t) & = b_3 \cos(x_1(t)) - b_4 x_3(t) + E + u(t)
	\end{align*}
	with constants $b_1 = 34.29$, $b_2 = 0.0$, $b_3 = 0.149$, $b_4 = 0.3341$, $P = 28.22$ and $E = 0.2405$, cf. \cite{GalazOrtegaBazanellaStankovic2003}. The equilibrium we wish to stabilize is located at $\xs = (1.124603730, 0, 0.9122974248)^\top$ and the running costs are defined as
	\begin{align*}
		\l(\x(t), \u(t)) = \| \x(t) - \xs \|_2^2 + \lambda \| u(t) \|_2^2
	\end{align*}
	where the parameter $\lambda = 0.01$ is used to penalize the taken control effort. Due to physical considerations $x_1$ and $x_3$ are restricted to the interval $[0, \pi/2]$ and $\R_{\geq 0}$ respectively.
\end{example}

%\karl{Kommasetzung, beachten wir sowas???} \juergen{Offensichtlich nicht...}
%
For Example \ref{Stability and Performance:ex:generator}, we want to determine an optimization horizon length $T$ such that asymptotic stability of the resulting MPC closed loop is guaranteed, i.e. $\overline{\alpha} = 0$. 

To achieve this goal, the following methodology can be applied to compute a function $B(\cdot)$ such that \eqref{Stability and Performance Bounds:thm:stability:eq1} holds:
\begin{enumerate}
	\item For each state $\x \in \bX$ derive a monotone function $B_{\x}: \R_{\geq 0} \rightarrow \R_{\geq 0}$ satisfying
		\begin{align}\label{NotationGrowthConditionLocal}
			V_t(\x) \leq B_{\x}(t) \l^\star(\x) \qquad\forall\, t \geq 0.
		\end{align}
	\item Then, $B(t)$ is defined pointwise as $\sup_{\x \in \bX} B_{\x}(t)$. 
\end{enumerate}
To accomplish Step (1) we exploit the fact that optimality of this bound is not needed. Hence, for each $t \in \N \cdot \Delta $ for some $\Delta > 0$ we solve
\begin{align*}
	\min_{\u \in \cU_{\x}([0,t))} J_t(\x,\u)
\end{align*}
over the class of sampled data systems with zero order hold and sampling period $\Delta$, i.e. control functions satisfying $u(t) = \text{\em constant}$ for $t \in [(n-1) \Delta, n \Delta)$, $n \in \N$. %Assuming admissibility of the resulting control function and dividing by $\l^\star(\x) = \| \x - \x^\star\|$ 
This yields monotone bounds $B_x(t)$, $t \in \Delta \N$, satisfying \eqref{NotationGrowthConditionLocal} by dividing by $\l^\star(\x) = \| \x - \x^\star\|$. Then, Step (2) is carried out in order to compute $B(t)$. Since neither terminal constraints nor costs were imposed, $V_t(\x)$ is monotonically increasing in $t$ which implies $V_t(\x) \leq V_{n \Delta}(\x) \leq B(n \Delta) \l^\star(\x)$, $t \in ((n-1)\Delta, n\Delta]$. Note that $B$ is only computed on a ``sufficiently'' large interval $[0,n^\star\Delta]$, $n^\star \in \N$, such that a desired performance estimate can be concluded by Theorem \ref{Stability and Performance Bounds:thm:stability}.
	
While the proposed procedure allows to verify Step (1) rigorously, we carry out Step (2) only approximately. The reason is twofold: Firstly, throughout this paper we assume feasibility for the set $\bX$, i.e. $\bX$ is control invariant, %$\bX$ is assumed to be a control invariant set, 
cf.~\cite{Blanchini1999} for a definition of control invariance. Hence, a suitable subset %set 
has to be computed a priori if this assumption is violated, see \cite[Chapter 8]{GruenePannek2011}. For Example \ref{Stability and Performance:ex:generator}, numerical computations show that an appropriately chosen level set $\cL$ %\cL \subset \bX$
of the value function $V_{0.6}$ %is contained in the interior of the state constrained set 
is control invariant. Secondly, this set has to be discretized. To this end, a grid of initial values $\mathcal{G}$ with discretization stepsize $\Delta_{x_i}$ within each direction $i \in \{1, \ldots, n\}$ is used. Then, $B(t)$ is approximately determined for the set $\cL$ by taking the supremum with respect to all states contained in the intersection $\overline{X} := \cL \cap \mathcal{G}$.

In conclusion, the proposed approach allows to rigorously ensure the relaxed Lyapunov inequality %for each state $x \in \cL$ 
except for the state space discretization assuming that a control invariant set %$\cL$ 
is given.
\begin{example}\label{Stability and Performance:ex:bound}
	Regarding Example \ref{Stability and Performance:ex:generator}, let $\Delta = 0.0125$ and the control horizon $\delta = 0.05$ be given. Additionally, we define $\mathcal{G} := [ \xs_1 - a_1, \xs_1 + a_1 ] \times [ \xs_2 - a_2, \xs_2 + a_2 ] \times [ \xs_3 - a_3, \xs_3 + a_3 ]$ with $a_1 = 0.4$, $a_2 = 0.5$, $a_3 = 0.9$ and discretization stepsize $\Delta_{x_i} := 0.02$, $i \in \{1, 2, 3\}$, and focus on the level set $\cL := \{ x \mid V_{0.6}(x) \leq 0.0081 \}$ whose convex hull satisfies all physical constraints, see Figure~\ref{Preliminaries:fig:set of initial values} for an illustration of $\overline{X} := \cL \cap \mathcal{G}$.
	\begin{figure}[!ht]
		\begin{center}
			\includegraphics[width=0.42\textwidth]{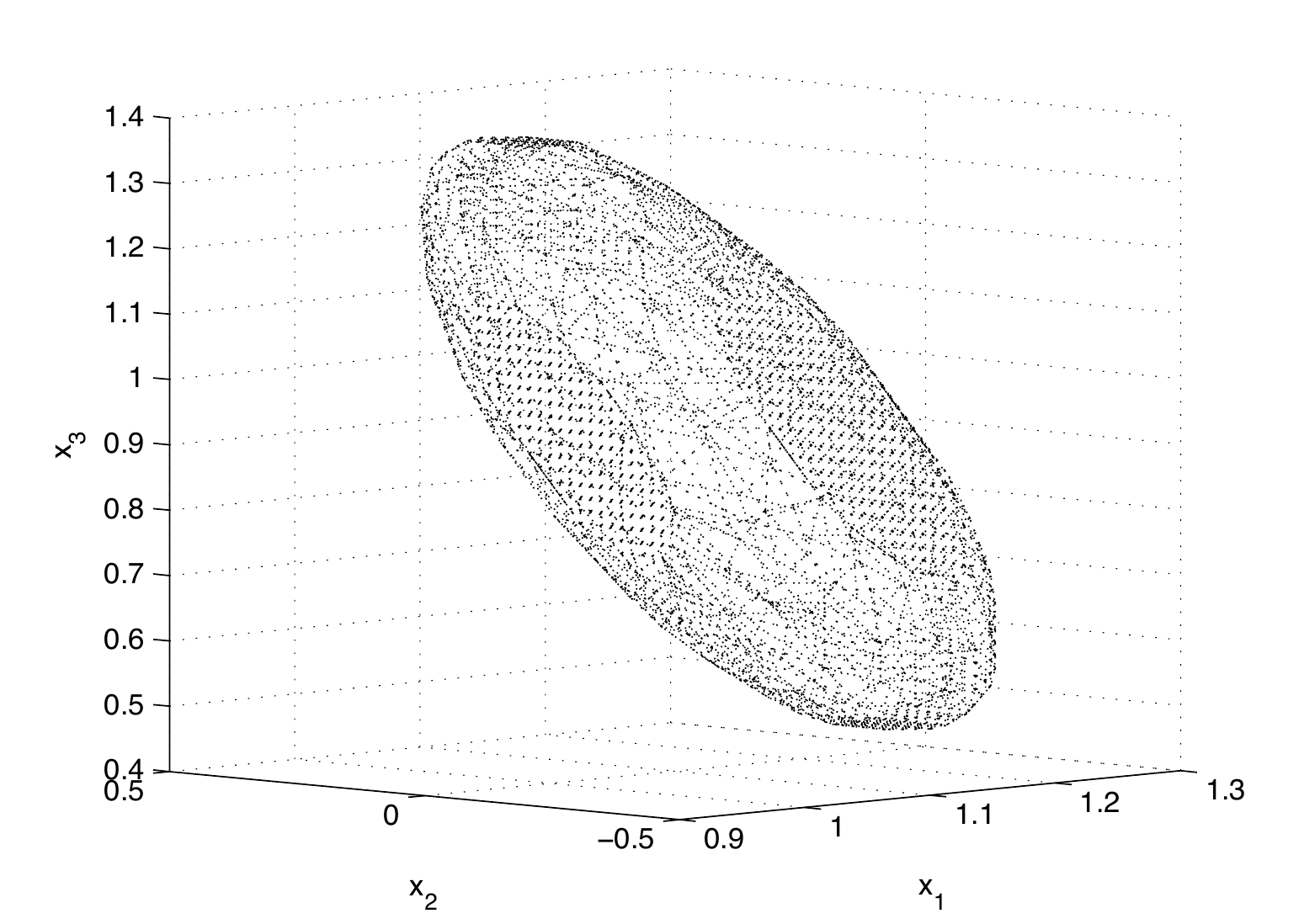}
			\caption{Illustration of the set $\overline{X} := \cL \cap \mathcal{G}$.}
			\label{Preliminaries:fig:set of initial values}
		\end{center}
	\end{figure}
	Then, Formula \eqref{Stability and Performance Bounds:thm:stability:eq2} enables us to determine an optimization horizon $T$ such that the relaxed Lyapunov inequality \eqref{Stability and Performance Bounds:thm:stability:eq3} holds with  $\alpha_{T,\delta} > \overline{\alpha} = 0$. For the computed function $B$ this methodology yields asymptotic stability of the MPC closed loop for optimization horizon $T = 2.6$.
\end{example}
%\karl{Bezug von $a_1$ und $a_2$ zu den Zustandsbeschr\"{a}nkungen: erw\"{a}hnen?}
\begin{remark}
	Note that the presented method verifies condition \eqref{Stability and Performance Bounds:thm:stability:eq1} for control functions $\u \in \mathcal{L}^\infty([0,T),\R^m)$, and allows to conclude asymptotic stability of the closed loop via Theorem \ref{Stability and Performance Bounds:thm:stability}, cf.~\cite[Remark 2.7]{WorthmannRebleGrueneAllgoewer2012}.
\end{remark}

% ------------------------------------------------------------------------
\section{Impact of the Control Horizon}
\label{Section:Increasing the Control Horizon}
% ------------------------------------------------------------------------

In Section \ref{Section:Stability and Performance Bounds} we showed how to ensure asymptotic stability for the proposed MPC scheme for a given $\delta$. In this section, we investigate %whether this quantity has an 
the impact of $\delta$ % on the deduced estimate 
on the required optimization horizon length $T$.

\begin{figure}[!ht]
	\begin{center}
		a) \includegraphics[width=0.2165\textwidth]{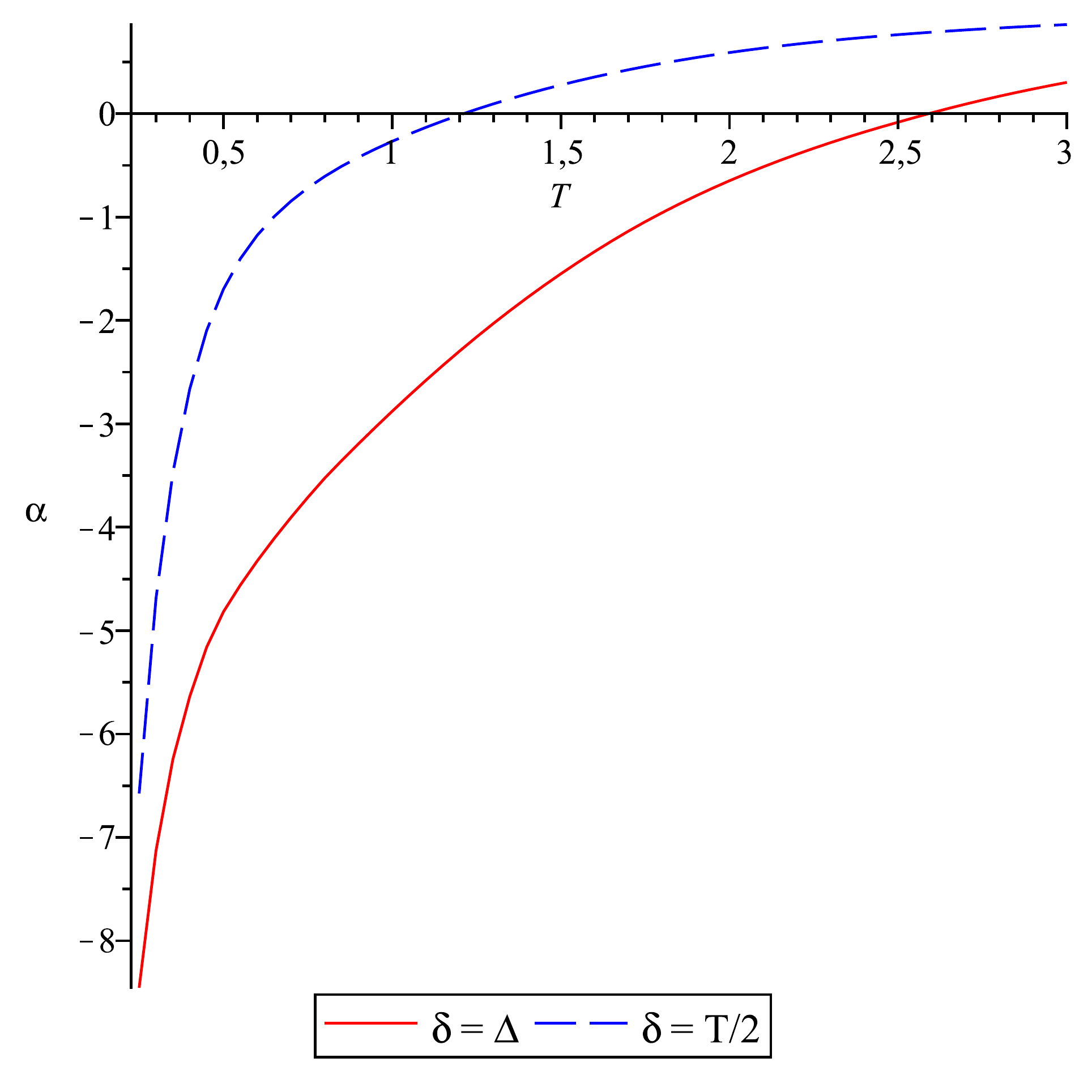} 
		b) \includegraphics[width=0.2165\textwidth]{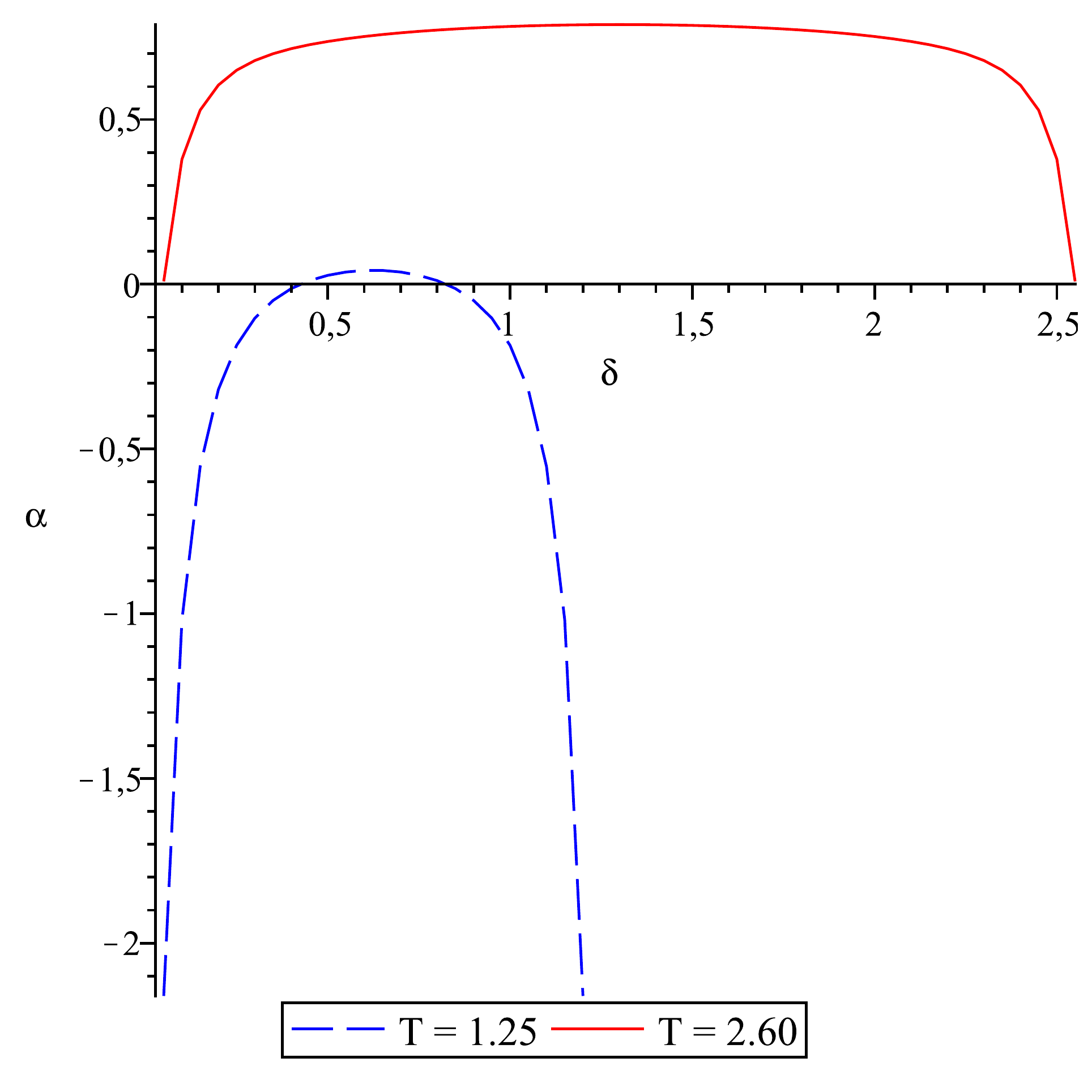}
		\caption{a) Development of $\alpha_{T, \delta}$ for varying $T$ and different choices of $\delta$. b) Development of $\alpha_{T, \delta}$ depending on the control horizon $\delta$.}
		\label{Increasing the Control Horizon:fig: horizons}
	\end{center}
\end{figure}

Considering Example \ref{Stability and Performance:ex:generator}, we compute $\alpha_{T,\Delta}$, $\Delta = 0.05$, and $\alpha_{T,T/2}$ for $T = n \Delta$, $n \in \{2,3,\ldots,60\}$, cf.~Figure~\ref{Increasing the Control Horizon:fig: horizons}a. Here, $\alpha_{T,\Delta} > 0$ holds for $T \geq 2.6 = 52 \Delta$ whereas this stability criterion holds for significantly shorter optimization horizons if the control horizon is chosen equal to $T/2$, that is $\alpha_{T,T/2} > 0$ for $T \geq 1.25 = 25 \Delta$. Hence, enlarging the control horizon seems to induce an improved performance index $\alpha_{T,\delta}$.

This numerical result motivates to investigate the influence of $\delta$ on $\alpha_{T,\delta}$. To this end, we fix the optimization horizon $T$ and compute $\alpha_{T,n\Delta}$ for $n \in \{1,2,\ldots,T/\Delta-1\}$, cf.~Figure~\ref{Increasing the Control Horizon:fig: horizons}b, which leads to the following observations: Firstly, a symmetry property seems to hold, i.e. $\alpha_{T,\delta} = \alpha_{T,T-\delta}$. Secondly, the performance estimates appear to increase up to the symmetry axis $T/2$. Both properties have been shown for systems which are exponentially controllable in terms of their stage costs, i.e. $\widetilde{B}(t) = C \int_0^t e^{-\mu s}\, ds$ for an overshoot constant $C \geq 1$ and a decay rate $\mu > 0$, cf.~\cite{GruenePannekWorthmann2012cdc}. Using the computed function $B$ instead of exponential decay we obtain $\widetilde{B}(t) \geq B(t)$ and therefore better horizon estimates as shown in \cite{Worthmann2011} for a discrete time setting. Despite the more general setting, symmetry still follows directly from Formula \eqref{Stability and Performance Bounds:thm:stability:eq2}. % can still be shown
\begin{corollary}\label{Increasing the Control Horizon:cor:symmetry}
	The performance estimate $\alpha_{T,\delta}$ given by Formula \eqref{Stability and Performance Bounds:thm:stability:eq2} satisfies $\alpha_{T, \delta} = \alpha_{T, T-\delta}$ for $\delta \in (0, T)$, i.e. $\alpha_{T,\delta}$ is symmetric with symmetry axis $\delta = T/2$.
\end{corollary}
%{\bf Proof}: Considering $\alpha_{T, \delta}$ and $\alpha_{T, T-\delta}$ from Formula \eqref{Stability and Performance Bounds:thm:stability:eq2} for $\delta \in (0, T)$ one can divide both fractions by their respective numerator showing the assertion. \qed

Unlike symmetry, we conjecture that there exists a counterexample negating monotonicity even if $B$ satisfies the additional condition \cite[Inequality (21)]{RebleAllgoewer2011} which is, however, violated in Example \ref{Stability and Performance:ex:generator}. For such an example \eqref{Stability and Performance Bounds:thm:stability:eq1} holds but $\alpha_{T, \delta}$ is not monotone on $(0,T/2]$, cf.~\cite{GruenePannekSeehaferWorthmann2010} for a counterexample in the discrete time setting.

Yet, for Example \ref{Stability and Performance:ex:generator} numerical results indicate that the monotonicity property holds and allows to conclude asymptotic stability of the MPC closed loop for significantly shorter optimization horizons.
\begin{example}\label{Increasing the Control Horizon:ex:generator2}
	Again consider Example \ref{Stability and Performance:ex:generator}. If we impose the horizon $T = 1.25$, then from Figure \ref{Increasing the Control Horizon:fig: horizons}b we observe that stability of the closed loop can only be guaranteed if $\delta \in [0.45, 0.8]$. If the horizon length is increased to $T = 2.60$, Figure~\ref{Increasing the Control Horizon:fig: horizons}a shows that for any chosen $\delta \in [0.05, 2.55]$ stability can be concluded.
\end{example}
In the context of arbitrary monotone and bounded functions $B: \R_{\geq 0} \rightarrow \R_{\geq 0}$ 
another interesting fact arises if we consider arbitrarily small $\delta$:%, \juergen{one can additionally show that asymptotic stability cannot be maintained for an arbitrarily small control horizon $\delta$:}
\begin{theorem}\label{Increasing the Control Horizon:thm:divergence}
	For any optimization horizon $T > 0$ we have that $\alpha_{T,\delta}$ %diverges to 
	goes to $- \infty$ for $\delta \rightarrow 0$. In particular, our stability condition $\alpha_{T,\delta} \geq 0$ cannot be maintained for arbitrarily small control horizon $\delta$.
\end{theorem}
{\bf Proof}: Follows directly from Formula \eqref{Stability and Performance Bounds:thm:stability:eq2}.\qed

Note that this assertion was solely shown for an exponentially controllable setting both for discrete and continuous time systems, cf.~\cite[Section 4]{RebleAllgoewer2011} and \cite[Section 5.1]{Worthmann2011}. Hence, the key contribution is the observation that this conclusion can also be drawn without the restriction to this particular class of systems.

Motivated by Example \ref{Increasing the Control Horizon:ex:generator2} and Theorem \ref{Increasing the Control Horizon:thm:divergence}, we propose an algorithm to obtain a control horizon length $\delta$ such that $\alpha_{T, \delta}$ exceeds a predefined suboptimality bound $\overline{\alpha}$. To this end, we introduce a partition $(\tau_k)_{k \in \{0, \ldots, n\}}$, $n \in \N_{> 1}$, of $[0, T]$ with $0 = \tau_0 < \tau_1 < \ldots < \tau_n = T$.  Such a setting naturally arises in the context of digital control for sampled data systems with zero order hold. Yet, we like to note that Algorithm \ref{Increasing the Control Horizon:alg:control horizon} is not limited to the digital control case. Additionally, we like to stress that monotonicity of $\alpha_{T, \delta}$ in $\delta$ is the center of this algorithm.

\begin{algorithm}
	\caption{MPC with increased control horizon}\label{Increasing the Control Horizon:alg:control horizon}
	Given: $T > 0$, $(\tau_k)_{k \in \{0, \ldots, n\}}$ with $n \in \N_{> 1}$ and $\overline{\alpha} \in [0, 1)$
	\begin{enumerate}
		\item Measure the current state $\hat{\x}$ 
		\item Set $k := 0$ and compute a minimizer $\us \in \cU_{\xh}([0, T))$ of \eqref{Preliminaries:eq:finite horizon cost functional} and $V_T(\xh) = J_T(\xh,\us)$.\\ Do
		\begin{enumerate}
			\item If $(k + 1) = n$: Set $\delta$ according to exit strategy and goto (3)
			\item Set $k := k + 1$, $\delta :=\tau_k$ and compute $V_T(\x_{\us}(\delta; \xh))$
			\item Compute $\alpha_{T, \delta}$, i.e.
				\begin{align}\label{Algorithm2ComputingAlpha}
					\alpha_{T, \delta} := \frac {V_T(\xh) - V_T(\x_{\us}(\delta; \xh))}{\int_0^\delta \hspace*{-0.1cm} \l(\x_{\us}(t;\xh), \us(t))\, dt}
				\end{align}
			%\item If $\alpha_{T, \delta} > \overline{\alpha}$: Goto (3)
		\end{enumerate}
		while $\alpha_{T, \delta} \leq \overline{\alpha}$
		\item Implement $\mu_{T, \delta}(t; \xh)|_{t \in [0, \delta)}$, shift the horizon forward in time by $\delta$ and goto (1)
	\end{enumerate}
\end{algorithm}

%The purpose of Algorithm \ref{Increasing the Control Horizon:alg:control horizon} \karl{hat ein Algorithmus einen ``purpose''} is to combine two aspects: 
Algorithm \ref{Increasing the Control Horizon:alg:control horizon} combines two aspects: the improved performance estimates obtained for larger control horizons $\delta$, and the inherent robustness resulting from using a feedback control law which benefits from updating the control law as often as possible. In order to illustrate this claim let us consider the numerical Example \ref{Stability and Performance:ex:generator} again. From Figure~\ref{Increasing the Control Horizon:fig: horizons}b we observed that asymptotic stability of the MPC closed loop can be shown for $T \geq 1.25$. Indeed, Algorithm \ref{Increasing the Control Horizon:alg:control horizon} only uses $\delta = \Delta = 0.05$ in each step independent of the chosen initial condition. Hence, MPC with constant control horizon $\delta = 0.05$ is performed safeguarded by our theoretically obtained estimates. Consequently, no exit strategy is needed since Step (2a) is excluded.

\begin{figure}[!ht]
	\begin{center}
		a) \includegraphics[width=0.2165\textwidth]{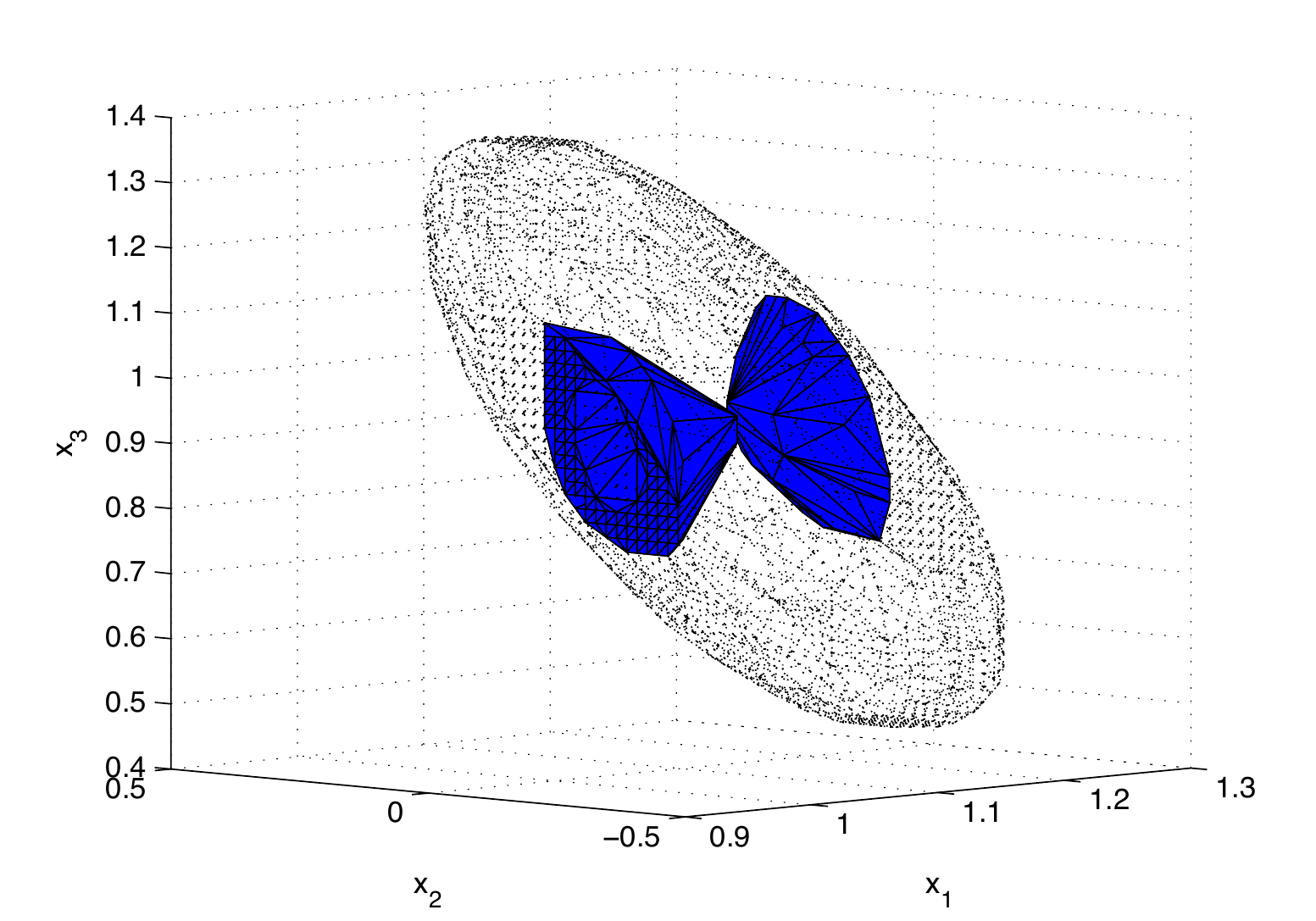}
		b) \includegraphics[width=0.2165\textwidth]{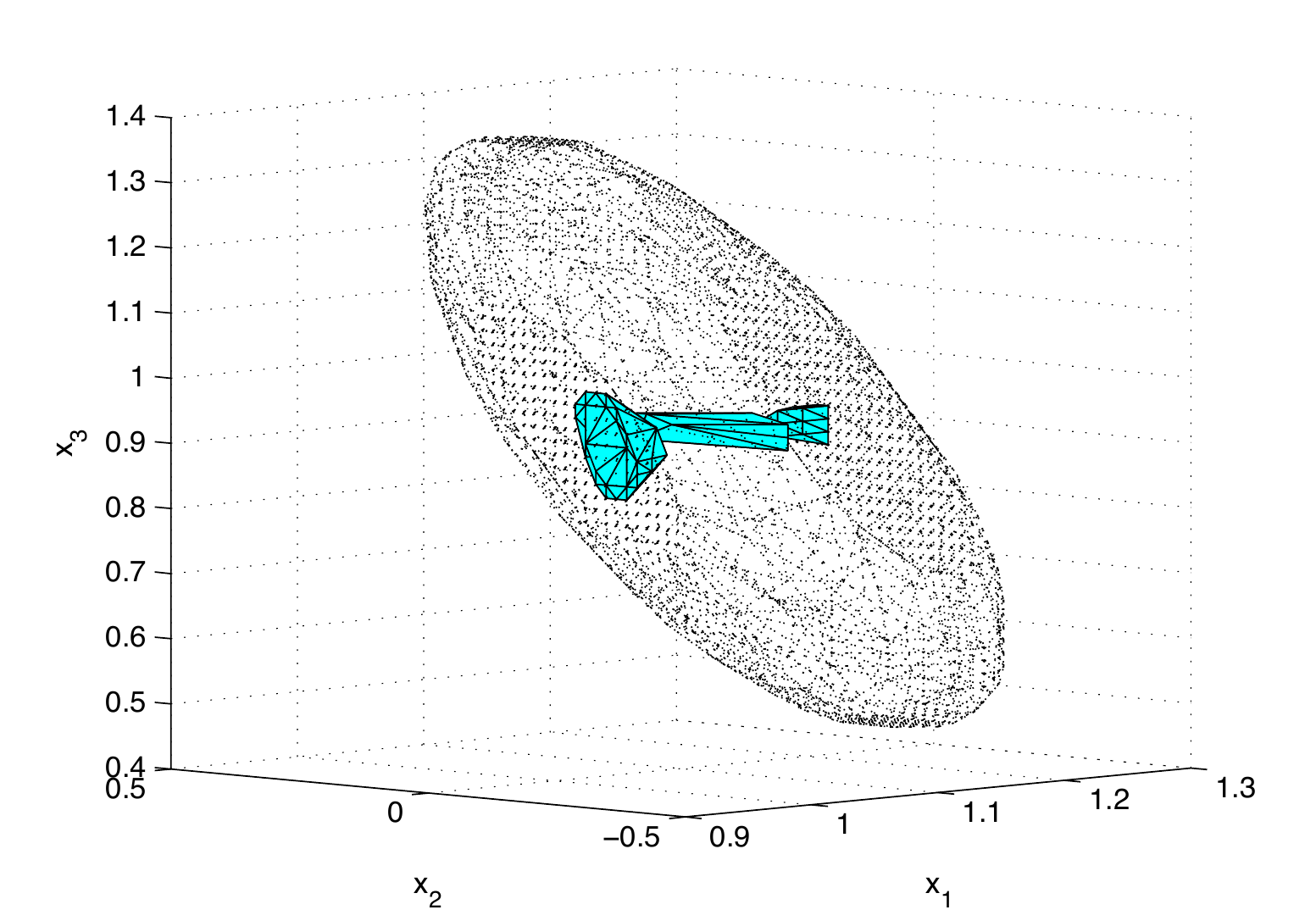}
		\caption{Set of initial values $x$ for which stability cannot be guaranteed for $T = 0.45$, $\Delta = 0.05$ with (a) $\delta = \Delta$ and (b) $\delta \in \{ \Delta, 2 \Delta \}$.}
		\label{Increasing the Control Horizon:fig:ex:N5}
	\end{center}
\end{figure}

We compute $\alpha_{T,\delta}$ by Equation \eqref{Algorithm2ComputingAlpha} for all $\x \in \overline{X}$ in order to test out the limits of Algorithm \ref{Increasing the Control Horizon:alg:control horizon}. %To verify that the observed benefits corresponding to using larger $\delta$ are not a theoretical relic, 
Here, $T = 0.45$ is the smallest optimization horizon such that $\alpha_{T,\Delta} = \alpha_{T,\Delta}(\x) > 0$ is satisfied for all $\x \in \overline{X}$. However, using Algorithm \ref{Increasing the Control Horizon:alg:control horizon} allows to ensure this conditions for $T = 0.25$ for varying control horizon $\delta \in \{\Delta,2\Delta,3\Delta,4\Delta\}$. Figure \ref{Increasing the Control Horizon:fig:ex:N5} shows the sets of states %vectors 
$\x \in \overline{X}$ for which $\alpha_{T, \delta}$ computed by \eqref{Algorithm2ComputingAlpha} is less than zero for the cases $\delta = \Delta$ and $\delta \in \{ \Delta, 2 \Delta\}$. Again, we see that if we allow for larger control horizons $\delta$, then the performance bound $\alpha_{T, \delta}$ increases, i.e. the set of state vectors for which stability cannot be guaranteed shrinks. Hence, for the considered Example \ref{Stability and Performance:ex:generator} Algorithm \ref{Increasing the Control Horizon:alg:control horizon} allows to drastically reduce the horizon $T$ while maintaining asymptotic stability of the closed loop.

The downside of considering potentially large control horizons $\delta$ is the possible lack of robustness in case of disturbances. Utilizing the introduced partition $(\tau_k)_{k \in \{0, \ldots, n\}}$, we can perform an update of the feedback law at time $\tau_j \in \{0, \ldots, n\}$ via
\begin{align}
	\label{Increasing the Control Horizon:eq:update control}
	\mu_{T, \delta}(t; x) := \u^\star(t - \tau_j; \x_{\mu_{T,\delta}}(\tau_j; x)), \quad t \geq \tau_j,
\end{align}
where we extended the notation of the open loop optimal control to $\u^\star(t; \x)$ to indicate which initial state is considered. Such an update can be applied whenever the following Lyapunov type update condition holds, see also \cite{PannekWorthmann2011a} for the discrete time setting.
\begin{proposition}\label{Increasing the Control Horizon:thm:update}
	Let $\overline{\alpha} \in [0, 1)$, $T > 0$ and a partition $(\tau_k)_{k \in \{0, \ldots, n\}}$, $n \in \N$, of $[0, T]$ with $0 = \tau_0 < \tau_1 < \ldots < \tau_n = T$ be given. For $\delta = \tau_k$ suppose \eqref{Stability and Performance Bounds:thm:stability:eq3} holds with $\alpha_{T, \delta} > \overline{\alpha}$ for some $\x \in \bX$. Let $\u^\star(\cdot; \x_{\mu_{T,\delta}}(\tau_j; \x))$ be a minimizer of \eqref{Preliminaries:eq:finite horizon cost functional} for some $j \in \{1, \ldots, k-1\}$. If additionally %condition
	\begin{align}
		& V_T(\x_{\u^\star}(\tau_k - \tau_j; \x_{\mu_{T,\delta}}(\tau_j; \x))) - V_{T - \tau_j}(\x_{\mu_{T,\delta}}(\tau_j; \x)) \nonumber \\
		& < ( 1 - \overline{\alpha} ) \int_0^{\tau_j} \l(\x_{\mu_{T,\delta}}(t; \x), \mu_{T,\delta}(t ; \x)) dt \label{Increasing the Control Horizon:thm:update:eq1} \\
		& - \overline{\alpha} \int_{0}^{\tau_k-\tau_j} \l(\x_{\u^\star}(t; \x_{\mu_{T,\delta}}(\tau_j; \x)), \u^\star(t; \x_{\mu_{T,\delta}}(\tau_j; \x))) dt \nonumber
	\end{align}
	holds% for some $\x \in \bX$
	, then the respective MPC feedback law $\mu_{T, \delta}(\cdot; \x)$ can be modified by \eqref{Increasing the Control Horizon:eq:update control} and the lower bound $\overline{\alpha}$ on the degree of suboptimality is locally maintained.
\end{proposition}
We like to stress that if the stabilizing partition index $k$ is known, then Proposition \ref{Increasing the Control Horizon:thm:update} allows for iterative updates of the feedback until this index is reached. Hence, only Step (3) of Algorithm \ref{Increasing the Control Horizon:alg:control horizon} needs to be adapted, see Algorithm \ref{Increasing the Control Horizon:alg:update} for a possible implementation.
\begin{algorithm}
	\caption{MPC with intermediate update}\label{Increasing the Control Horizon:alg:update}
	\begin{enumerate}
		\setcounter{enumi}{2}
		\item For $n = 0, \ldots k-1$ do
		\begin{enumerate}
			\item Implement $\mu_{T, \delta}(t; \xh)|_{t \in [\tau_n, \tau_{n+1})}$
			\item Compute $\u^\star(\cdot; \x_{\mu_{T,\delta}}(\tau_{n+1}; \hat{\x}))$ 
			\item If \eqref{Increasing the Control Horizon:thm:update:eq1} holds: Update $\mu_{T, \delta}(\cdot; \xh)$ via \eqref{Increasing the Control Horizon:eq:update control} 
		\end{enumerate}
		Shift the horizon forward in time by $\delta$ and goto (1)
	\end{enumerate}
\end{algorithm}

%\karl{Bitte \"{U}berschriften checken!}

% ------------------------------------------------------------------------
\section{Aggregated Performance}
\label{Section:Slack Stability}
% ------------------------------------------------------------------------

Usually, the relaxed Lyapunov inequality \eqref{Stability and Performance Bounds:thm:stability:eq3} is tight for only a few points in the state space. Hence, if a closed loop trajectory $\x_{\mu_{T, \delta}}(\cdot; \x)$ visits a point for which \eqref{Stability and Performance Bounds:thm:stability:eq3} is not an equality, then we can compute the occuring slack along the closed loop via
\begin{align}
	s(t, \x) := & V_T(\x) - V_T(\x_{\mu_{T, \delta}}(t; \x)) \nonumber \\
	\label{Slack Stability:eq:slack}
	& - \overline{\alpha} \int_0^{t} \l(\x_{\mu_{T, \delta}}(s; \x), \u_{\text{MPC}}^{\mu_{T, \delta}}(s))\, ds.
\end{align}
This slack can be used to weaken the requirement \eqref{Stability and Performance Bounds:thm:stability:eq3} by considering the closed loop instead of the open loop solution. For simplicity of exposition, we formulate the following result using constant $\delta$, yet the conclusion also holds in the context of time varying control horizons as in Algorithm \ref{Increasing the Control Horizon:alg:control horizon}.
\begin{theorem}\label{Slack Stability:thm:convergence}
	Consider an admissible feedback law $\u_{\text{MPC}}^{\mu_{T, \delta}}$, an initial value $\x_{0} \in \bX$ and $\overline{\alpha} \in (0, 1)$ to be given. Furthermore, suppose there exist $\K_\infty$ functions $\underline{\eta}$, $\overline{\eta}$ such that $\l^\star(x) \geq \underline{\eta}(\| \x - \x^\star \|)$ and $V_T(\x) \leq \overline{\eta}(\| \x - \x^\star \|)$ hold for all $\x \in \bX$. If additionally $s(t, \x)$ from \eqref{Slack Stability:eq:slack} converges for $t$ tending to infinity, then the MPC closed loop trajectory with initial value $\x_0$ behaves like an asymptotically stable solution. Furthermore, the following performance estimate holds
	\begin{align}\label{InequalitySlackTheorem}
		V_\infty^{\mu_{T, \delta}}(\x_0) \leq \frac{1-\lim_{t \rightarrow \infty}s(t,\x_0)/V_T(\x_0)}{\overline{\alpha}} V_\infty(\x_0).
	\end{align}
\end{theorem}
{\bf Proof}: Let the limit of $s(t,\x_0)$ for $t \rightarrow \infty$ be denoted by $\theta$. Then, for given $\varepsilon > 0$, there exists a time instant $t^\star$ such that $\| s(t,\x_0) - \theta \| \leq \varepsilon$ holds for all $t \geq t^\star$. As a consequence, we obtain
\begin{align*}
	s(t,\x_0) = & s(t^\star,\x_0) - V_T(\x_{\mu_{T,\delta}}(t;\x_0)) + V_T(\x_{\mu_{T,\delta}}(t^\star;\x_0)) \\
	& - \overline{\alpha} \int_{t^\star}^t \l(\x_{\mu_{T,\delta}}(s;\x_0),u_{\text{MPC}}^{\mu_{T,\delta}}(s))\, ds \geq \theta - \varepsilon.
\end{align*}
Since $V_T(\x_{\mu_{T,\delta}}(t^\star;\x_0)) \leq \overline{\eta}(\| \x_{\mu_{T,\delta}}(t^\star;\x_0) - \x^\star \|)$ as well as $V_T(\x_{\mu_{T,\delta}}(t;\x_0)) \geq 0$ hold, this inequality implies
\begin{align*}
	& 2 \varepsilon + \overline{\eta}(\| \x_{\mu_{T,\delta}}(t^\star;\x_0) - \x^\star \|) \\
	\geq & s(t^\star,\x_0) - \theta + \varepsilon + V_T(\x_{\mu_{T,\delta}}(t^\star;\x_0)) - V_T(\x_{\mu_{T,\delta}}(t;\x_0)) \\
	\geq & \overline{\alpha} \int_{t^\star}^t \l(\x_{\mu_{T,\delta}}(s;\x_0),u_{\text{MPC}}^{\mu_{T,\delta}}(s))\, ds.
\end{align*}
Hence, boundedness of the integral on the right hand side can be concluded. Now, due to positivity and continuity of $\l$ we have $\lim_{t \rightarrow \infty} \l(\x_{\mu_{T,\delta}}(t;\x_0),u_{\text{MPC}}^{\mu_{T,\delta}}(t)) = 0$. In turn, the latter ensures $\underline{\eta}(\| \x_{\mu_{T,\delta}}(t;\x_0) - \x^\star \|) \rightarrow 0$ and, thus, $\x_{\mu_{T,\delta}}(t;\x_0) \rightarrow \x^\star$ for $t$ approaching infinity. Inequality \eqref{InequalitySlackTheorem} is shown by using
\begin{align*}
	\lim_{t \rightarrow \infty} V_T(\x_{\mu_{T, \delta}}(t; \x)) \leq \lim_{t \rightarrow \infty} \overline{\eta}(\|\x_{\mu_{T, \delta}}(t; \x) - \x^\star\|) = 0
\end{align*}
in combination with \eqref{Slack Stability:eq:slack} to obtain
\begin{eqnarray*}
	\overline{\alpha} V_\infty^{\mu_{T, \delta}}(\x) & = & \lim_{t \rightarrow \infty} \overline{\alpha} \int_0^{t} \l(\x_{\mu_{T, \delta}}(s; \x), \u_{\text{MPC}}^{\mu_{T, \delta}}(s))\, ds \\
	& = & V_T(\x) - \lim_{t \rightarrow \infty} s(t,\x).
\end{eqnarray*}
Then, since $V_T(\x) \leq V_\infty(\x)$, the assertion follows.\qed

Note that within Theorem \ref{Slack Stability:thm:convergence} we did not assume semipositivity but convergence of $s(t, \x)$ to conclude stability. Here, we like to stress that $\lim_{t \rightarrow \infty} s(t,\x) \geq 0$ implies a suboptimality index $\alpha_{T, \delta} \geq \overline{\alpha}$. Clearly, both the stability and performance result shown in Theorem \ref{Slack Stability:thm:convergence} can be extended to assertions for all $\x \in \bX$ if $s(t,\x)$ converges for every $\x \in \bX$ or a uniform lower bound can be found, i.e. $\inf_{\x \in \bX} \lim_{t \rightarrow \infty} s(t,\x) = \theta > - \infty$.

Apart from its theoretical impact, $s(t,\x)$ is also meaningful at runtime of the MPC algorithm. For instance, the condition $s(t, \x) \geq 0$ can be checked at each time instant $t = n\delta$, $n \in \N$, instead of $\alpha_{T,\delta} > \overline{\alpha}$, cf.~\eqref{Algorithm2ComputingAlpha}. This is particularly useful since accumulated slack can be used in order to compensate local violations of $\alpha_{T,\delta} \geq \overline{\alpha}$, i.e. weakening the stability condition \eqref{Algorithm2ComputingAlpha}, as long as the overall performance is still satisfactory.

\begin{figure}[!ht]
	\begin{center}
		\includegraphics[width=0.42\textwidth]{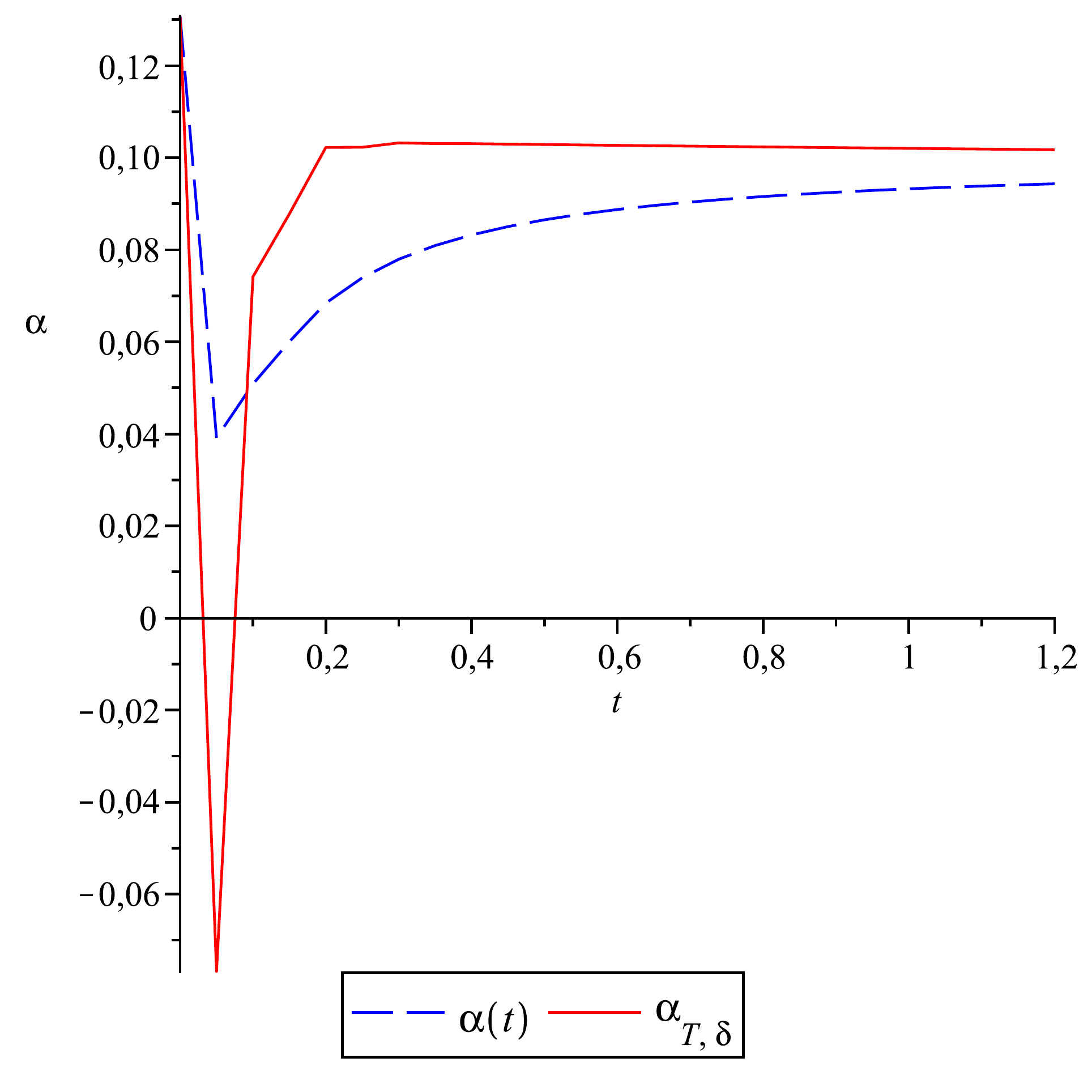}
		\caption{Development of $\alpha_{T, \delta}$ and $\alpha(t)$ for $T = 0.25$ and $\delta = 0.05$ along a specific closed loop solution.}
		\label{Slack Stability:fig:alpha}
	\end{center}
\end{figure}

If $s(t, \x) < 0$ occurs within the MPC algorithm, the slack can also be used to form an exit strategy. To this end, we denote the performance of the MPC closed loop until time $t$ by
\begin{align}\label{NotationCurrentPerformance}
	\alpha(t) := \frac {V_T(\x) - V_T(\x_{\mu_{T, \delta}}(t; \x))}{\int_0^{t} \l(\x_{\mu_{T, \delta}}(s; \x), \u_{\text{MPC}}^{\mu_{T, \delta}}(s))\, ds}.
\end{align}
Now, if $s(t,\x) < 0$ but 
$\alpha(t) > 0$, then stability is still maintained, yet the current performance index is worse than the desired bound $\overline{\alpha}$.
\begin{example}
	Again we consider Example \ref{Stability and Performance:ex:generator} and analyze the MPC closed loop for $\mu_{T, \delta}$ with $T = 0.25$ and $\delta = \Delta = 0.05$. If we choose the initial value $\x_0 = (1.03960373, -0.085, 0.9122974248)^\top$, then we observe from Figure \ref{Slack Stability:fig:alpha} that the local estimate $\alpha_{T, \delta}$ drops below zero for $t = \delta$, i.e. stability cannot be  guaranteed for $\delta = \Delta$. Yet, computing $\alpha(t)$ according to \eqref{NotationCurrentPerformance} shows that the relaxed Lyapunov inequality is satisfied after two steps of the MPC algorithm. Hence, using the slack information %stored in $s(t,\x)$ and 
	incorporated in $\alpha(t)$ allows to conclude asymptotic stability.
\end{example}

% ------------------------------------------------------------------------
\section{Conclusions and Outlook}
\label{Section:Conclusions and Outlook}
% ------------------------------------------------------------------------
We have shown a methodology to verify the assumptions introduced in \cite{RebleAllgoewer2011} under which stability of the 
%unconstrained 
MPC closed loop without terminal constraints can be guaranteed. Additionally we presented an algorithmic approach %extension 
for varying control horizons which allows %allowed us 
to reduce the optimization horizon length. Last, we provided robustification methods via update and slack based rules.

% ------------------------------------------------------------------------
\begin{ack}                              
	This work was supported by the DFG priority research program 1305 ``Control Theory of Digitally Networked Dynamical Systems'', grant. no. Gr1569/12-1, and by the German Federal Ministry of Education and Research (BMBF) project ``SNiMoRed'', grant no. 03MS633G.
\end{ack}
% ------------------------------------------------------------------------

%% ------------------------------------------------------------------------
%\bibliographystyle{ifacconf-harvard}
%\bibliography{ifacconf}

\begin{thebibliography}{22}
\providecommand{\natexlab}[1]{#1}
\providecommand{\url}[1]{\texttt{#1}}
\expandafter\ifx\csname urlstyle\endcsname\relax
  \providecommand{\doi}[1]{doi: #1}\else
  \providecommand{\doi}{doi: \begingroup \urlstyle{rm}\Url}\fi

\bibitem[Alamir and Bornard(1995)]{AlamirBornard1995}
M.~Alamir and G.~Bornard.
\newblock Stability of a truncated infinite constrained receding horizon
  scheme: the general discrete nonlinear case.
\newblock \emph{Automatica}, 31\penalty0 (9):\penalty0 1353--1356, 1995.

\bibitem[Bardi and Capuzzo-Dolcetta(1997)]{BardiCapuzzo1997}
M.~Bardi and I.~Capuzzo-Dolcetta.
\newblock \emph{Optimal control and viscosity solutions of
  {H}amilton-{J}acobi-{B}ellman equations}.
\newblock Birkh\"auser, 1997.

\bibitem[Blanchini(1999)]{Blanchini1999}
F.~Blanchini.
\newblock {S}et {I}nvariance in {C}ontrol.
\newblock \emph{Automatica}, 35:\penalty0 1747--1767, 1999.

\bibitem[Camacho and Bordons(2004)]{CamachoBordons2004}
E.~F. Camacho and C.~Bordons.
\newblock \emph{{Model Predictive Control}}.
\newblock Springer, 2004.

\bibitem[Chen and Allg{\"o}wer(1998)]{ChenAllgoewer1998}
H.~Chen and F.~Allg{\"o}wer.
\newblock A quasi-infinite horizon nonlinear model predictive control scheme
  with guaranteed stability.
\newblock \emph{Automatica}, 34\penalty0 (10):\penalty0 1205--1218, 1998.

\bibitem[Galaz et~al.(2003)Galaz, Ortega, Bazanella, and
  Stankovic]{GalazOrtegaBazanellaStankovic2003}
M.~Galaz, R.~Ortega, A.~S. Bazanella, and A.~M. Stankovic.
\newblock An energy-shaping approach to the design of excitation control of
  synchronous generators.
\newblock \emph{Automatica}, 39\penalty0 (1):\penalty0 111--119, 2003.

\bibitem[Giselsson(2010)]{Giselsson2010}
P.~Giselsson.
\newblock {A}daptive {N}onlinear {M}odel {P}redictive {C}ontrol with
  {S}uboptimality and {S}tability {G}uarantees.
\newblock In \emph{Proceedings of the 49th Conference on Decision and Control,
  Atlanta, GA, USA}, pages 3644--3649, 2010.

\bibitem[Gr{\"u}ne(2009)]{Gruene2009}
L.~Gr{\"u}ne.
\newblock Analysis and design of unconstrained nonlinear {MPC} schemes for
  finite and infinite dimensional systems.
\newblock \emph{SIAM Journal on Control and Optimization}, 48:\penalty0
  1206--1228, 2009.

\bibitem[Gr\"une and Pannek(2011)]{GruenePannek2011}
L.~Gr\"une and J.~Pannek.
\newblock \emph{Nonlinear Model Predictive Control: Theory and Algorithms}.
\newblock Springer, 2011.

\bibitem[Gr\"une et~al.(2010)Gr\"une, Pannek, Seehafer, and
  Worthmann]{GruenePannekSeehaferWorthmann2010}
L.~Gr\"une, J.~Pannek, M.~Seehafer, and K.~Worthmann.
\newblock Analysis of unconstrained nonlinear {M}{P}{C} schemes with varying
  control horizon.
\newblock \emph{SIAM Journal on Control and Optimization}, 48\penalty0
  (8):\penalty0 4938--4962, 2010.

\bibitem[Gr\"une et~al.(2012)Gr\"une, Pannek, and
  Worthmann]{GruenePannekWorthmann2012cdc}
L.~Gr\"une, J.~Pannek, and K.~Worthmann.
\newblock {Ensuring Stability in Networked Systems with Nonlinear MPC for
  Continuous Time Systems}.
\newblock In \emph{Proceedings of the 51st Converence on Decision and Control,
  Maui, Hawaii, USA}, 2012.
\newblock To appear.

\bibitem[Jadbabaie and Hauser(2005)]{JadbabaieHauser2005}
A.~Jadbabaie and J.~Hauser.
\newblock On the stability of receding horizon control with a general terminal
  cost.
\newblock \emph{IEEE Trans. Automat. Control}, 50\penalty0 (5):\penalty0
  674--678, 2005.
\newblock ISSN 0018-9286.

\bibitem[Keerthi and Gilbert(1988)]{KeerthiGilbert1988}
S.~S. Keerthi and E.~G. Gilbert.
\newblock Optimal infinite horizon feedback laws for a general class of
  constrained discrete-time systems: stability and moving horizon
  approximations.
\newblock \emph{Journal of Optimization Theory and Applications}, 57:\penalty0
  265--293, 1988.

\bibitem[Maciejowski(2002)]{Maciejowski2002}
J.~M. Maciejowski.
\newblock \emph{{P}redictive {C}ontrol with {C}onstraints}.
\newblock Prentice-Hall, Harlow, England, 2002.

\bibitem[Pannek and Worthmann(2011)]{PannekWorthmann2011a}
J.~Pannek and K.~Worthmann.
\newblock {MPC Algorithms with Stability and Performance Guarantees}.
\newblock http://arxiv.org/pdf/1109.6153.pdf, 2011.

\bibitem[Primbs and Nevisti\'c(2000)]{PrimbsNevistic2000}
J.A. Primbs and V.~Nevisti\'c.
\newblock Feasibility and stability of constrained finite receding horizon
  control.
\newblock \emph{Automatica}, 36:\penalty0 965--971, 2000.

\bibitem[Rawlings and Mayne(2009)]{RawlingsMayne2009}
J.~B. Rawlings and D.~Q. Mayne.
\newblock \emph{Model Predictive Control: Theory and Design}.
\newblock Nob Hill Publishing, 2009.

\bibitem[Reble and Allg\"{o}wer(2011)]{RebleAllgoewer2011}
M.~Reble and F.~Allg\"{o}wer.
\newblock {Unconstrained Model Predictive Control and Suboptimality Estimates
  for Nonlinear Continuous-Time Systems}.
\newblock \emph{Automatica}, 2011.
\newblock Accepted.

\bibitem[Tuna et~al.(2006)Tuna, Messina, and Teel]{TunaMessinaTeel2006}
S.E. Tuna, M.J. Messina, and A.R. Teel.
\newblock {S}horter horizons for model predictive control.
\newblock In \emph{Proceedings of the 2006 American Control Conference,
  Minneapolis, Minnesota, USA}, 2006.

\bibitem[Worthmann(2011)]{Worthmann2011}
K.~Worthmann.
\newblock \emph{{S}tability {A}nalysis of {U}nconstrained {R}eceding {H}orizon
  {C}ontrol {S}chemes}.
\newblock PhD thesis, University of Bayreuth, 2011.

\bibitem[Worthmann(2012)]{Worthmann2012nmpc}
K.~Worthmann.
\newblock {Estimates on the Required Prediction Horizon in MPC: A Numerical
  Case Study}.
\newblock In \emph{Proceedings of the Conference on Nonlinear Model Predicitve
  Control, Nordwijkerhout, the Netherlands}, 2012.
\newblock To appear.

\bibitem[Worthmann et~al.(2012)Worthmann, Reble, Gr\"une, and
  Allg{\"o}wer]{WorthmannRebleGrueneAllgoewer2012}
K.~Worthmann, M.~Reble, L.~Gr\"une, and F.~Allg{\"o}wer.
\newblock The role of sampling for stability and performance in unconstrained
  model predictive control.
\newblock 2012.
\newblock Preprint, University of Bayreuth. Submitted.

\end{thebibliography}
%% ------------------------------------------------------------------------

%% ------------------------------------------------------------------------
%\appendix
%% ------------------------------------------------------------------------

% ------------------------------------------------------------------------
\end{document}